\documentclass[11pt,reqno]{amsart}
\usepackage{amssymb,amsmath,amsthm,newlfont,enumerate}

\theoremstyle{plain}
\newtheorem{theorem}{Theorem}[section]
\newtheorem{proposition}[theorem]{Proposition}
\newtheorem{lemma}[theorem]{Lemma}
\newtheorem{corollary}[theorem]{Corollary}

\theoremstyle{definition}

\theoremstyle{remark}

\newcommand{\CC}{\mathbb{C}}

\newcommand{\FF}{\mathbb{F}}
\newcommand{\RR}{\mathbb{R}}
\newcommand{\ZZ}{\mathbb{Z}}

\renewcommand{\tilde}{\widetilde}
\DeclareMathOperator{\supp}{supp}

\begin{document}

\title[Weighted convolution on the integers]{A counterexample
concerning weighted convolution on the integers}

\author[Morneau-Gu\'erin]{Fr\'ed\'eric Morneau-Gu\'erin}
\address{D\'epartement  \'Education, Universit\'e T\'ELUQ, 455 rue du Parvis, Qu\'ebec (QC), G1K 9H6, Canada}
\email{frederic.morneau-guerin@teluq.ca}

\author[Ransford]{Thomas Ransford}
\address{D\'epartement de math\'ematiques et de statistique, Universit\'e Laval, Qu\'ebec (QC), G1V 0A6, Canada}
\email{ransford@mat.ulaval.ca}

\date{2 Nov 2026}

\begin{abstract}
Let $1<p<\infty$ and let $w:\ZZ\to(0,\infty)$ be a weight.
When is the weighted space $\ell^p(\ZZ,w)$ an algebra with respect to convolution? A well-known sufficient condition is that 
$(w^{-q}*w^{-q})/w^{-q}$ be bounded on $\ZZ$,
where $q:=p/(p-1)$ is the index conjugate to $p$.
We present an example to show that this sufficient
condition is not necessary.
 \end{abstract}

\thanks{Research of second author  supported  by NSERC Discovery Grant RGPIN-2026-04565}


\makeatletter
\@namedef{subjclassname@2020}{\textup{2020} Mathematics Subject Classification}
\makeatother


\maketitle

\section{Introduction}\label{S:Intro}

Let $1\le p<\infty$ and let $w:\ZZ\to(0,\infty)$ be a weight, i.e., a strictly positive function.
The corresponding weighted $\ell^p$-space is
\[
\ell^p(\ZZ,w):=\Bigl\{f:\ZZ\to\CC:~\|f\|_{p,w}:=\Bigl(\sum_{x\in\ZZ}|f(x)|^pw(x)^p\Bigr)^{1/p}<\infty\Bigr\}.
\]

As explained in  \cite{Mo26} and \cite{MR20},
for various applications in harmonic and functional analysis,
it is important to know
for which weights $\ell^p(\ZZ,w)$ is stable with
respect to convolution. By this, we mean that, whenever
$f,g\in\ell^p(\ZZ,w)$,  their convolution product
\[
(f*g)(x):=\sum_{y\in\ZZ}f(y)g(x-y)
\]
converges absolutely for all $x\in\ZZ$ and satisfies $f*g\in\ell^p(\ZZ,w)$. If this is the case, then a standard
functional-analysis argument
(see, for example, \cite[Proposition~3.3]{MR20})
shows that the norm $\|\cdot\|_{p,w}$
is submultiplicative, in the sense that
\[
A_p(w):=\sup\{\|f*g\|_{p,w}:~\|f\|_{p,w}\le1,\,\|g\|_{p,w}\le1\}<\infty,
\] 
and therefore that $\ell^p(\ZZ,w)$ is a Banach algebra.

It is well known and easy to see that, if $p=1$, then a necessary and sufficient condition
for $\ell^1(\ZZ,w)$ to be an algebra with respect to 
convolution is that 
\[
B_1(w):=\sup_{x,y\in\ZZ}\frac{w(x)}{w(y)w(x-y)}<\infty,
\]
and in this case we have $A_1(w)\le B_1(w)$.

For $p>1$, the corresponding result is the following theorem.

\begin{theorem}\label{T:basic}
Let $1<p<\infty$ and let $w:\ZZ\to(0,\infty)$ be a weight. If
$(w^{-q}*w^{-q})/w^{-q}$ is bounded on $\ZZ$,
where $q:=p/(p-1)$ is the conjugate index to $p$, 
then $\ell^p(\ZZ,w)$ is an algebra with respect
to convolution. Moreover, we have $A_p(w)\le B_p(w)$, where
\[
B_p(w):=\Bigl\|\frac{w^{-q}*w^{-q}}{w^{-q}}\Bigr\|_\infty^{1/q}=
\Bigl(\sup_{x\in\ZZ}\sum_{y\in\ZZ}\Bigl(\frac{w(x)}{w(y)w(x-y)}\Bigr)^q\Bigr)^{1/q}.
\]
\end{theorem}

This result is also well known; a simple proof can be found,
for example, in \cite[Theorem~1.1]{MR20}.
However, it differs from the case $p=1$ in that it is merely
a \emph{sufficient} condition for $\ell^p(\ZZ,w)$ to be an algebra.
Until now, it was not known whether it is also \emph{necessary}.
Our goal in this article is to establish the following result showing that this sufficient condition is not necessary.

\begin{theorem}\label{T:main}
For each $1<p<\infty$, there exists a weight $w$ on $\ZZ$
such that $\ell^p(\ZZ,w)$ is an algebra with respect to convolution, 
but
\[
\sup_{x\in\ZZ}\sum_{y\in\ZZ}\Bigl(\frac{w(x)}{w(y)w(x-y)}\Bigr)^q=\infty.
\]
\end{theorem}

The novelty of this theorem is that the underlying group is~$\ZZ$. 
Analogous results had previously been obtained for $\ell^p(\ZZ^+,w)$ (see \cite{Fr78}), for  $L^p(\RR,w)$ (see \cite{Ku06}),  and for $\ell^p(\FF_2,w)$, where $\FF_2$ is the free group on two generators (see \cite{Ku09}).
However, the counterexamples in each of these cases made
use of specific properties of the domain not shared by $\ZZ$:
that $\ZZ^+$ is a semigroup and not a group; that $\RR$ is not a discrete group; and that $\FF_2$ is not an abelian group.
Because of this,  none of the examples could be extended to cover the case of~$\ZZ$, which had remained open until now.

More recently, it was shown in \cite{Mo26} that there exists a discrete abelian group $G$ and a weight $w$ on $G$ such that
$\ell^2(G,w)$ is an algebra, but 
\[
\sup_{x\in G}\sum_{y\in G}\Bigl(\frac{w(x)}{w(y)w(x-y)}\Bigr)^2=\infty.
\]
The construction depends on the fact that the group in question is a torsion group (every element has finite order), so it is not immediately applicable to the case of~$\ZZ$. Nonetheless, it turns out that the ideas in \cite{Mo26} can be used as a basis for a proof of Theorem~\ref{T:main}. 
This is the route that we shall follow.

We conclude this introduction by citing two well-known inequalities 
that will be used without comment in what follows.
First, for pointwise products, we have H\"older's inequality: if $1<p<\infty$ and $q=p/(p-1)$, if $f\in\ell^p(\ZZ)$ and if $g\in\ell^q(\ZZ)$, then
$fg\in\ell^1(\ZZ)$ with
\begin{equation}\label{E:Holder}
\|fg\|_1\le \|f\|_p\|g\|_q.
\end{equation}
Second, for convolution products, we have (a special case of) Young's inequality: if
$1\le p<\infty$, if $f\in\ell^1(\ZZ)$ and if $g\in \ell^p(\ZZ)$,
then $f*g\in\ell^p(\ZZ)$ with
\begin{equation}\label{E:Young}
\|f*g\|_p\le \|f\|_1\|g\|_p.
\end{equation}

%%%%%%%%%%%%%%%%%%%%%%%%%%%%%%%%%

\section{Proof strategy}\label{S:strategy}
Throughout this section, we 
fix $1<p<\infty$, and set 
$q:=p/(p-1)$ and $s:=1/q+2/3$.
We write $w_0$ for the weight on $\ZZ$
defined by
$w_0(x):=(1+|x|)^s$.

\subsection{Preliminary remarks}\label{S:prelim}
If $x,y\in\ZZ$, then
\[
\frac{1+|x|}{(1+|y|)(1+|x-y|)}
\le \frac{1+|y|+|x-y|}{(1+|y|)(1+|x-y|)}
\le \frac{1}{1+|y|}+\frac{1}{1+|x-y|}.
\]
Hence, for all $x\in\ZZ$,
\begin{align*}
\sum_{y\in\ZZ}\Bigl(\frac{w_0(x)}{w_0(y)w_0(x-y)}\Bigr)^q
&=\sum_{y\in\ZZ}\Bigl(\frac{1+|x|}{(1+|y|)(1+|x-y|)}\Bigr)^{qs}\\
&\le \sum_{y\in\ZZ}\Bigl( \frac{1}{1+|y|}+\frac{1}{1+|x-y|}\Bigr)^{qs}\\
&\le 2^{qs-1}\sum_{y\in\ZZ}\Bigl( \frac{1}{1+|y|}\Bigr)^{qs}+
2^{qs-1}\sum_{y\in\ZZ}\Bigl(\frac{1}{1+|x-y|}\Bigr)^{qs}\\
&=2^{qs}\sum_{y\in\ZZ}\Bigl( \frac{1}{1+|y|}\Bigr)^{qs}<\infty.
\end{align*}
It follows that $B_p(w_0)<\infty$.
More generally, if $w$ is any weight on $\ZZ$ such that
$w=w_0$ outside a finite subset of $\ZZ$, then, 
since $w/w_0$ is bounded
above and below, we have $B_p(w)<\infty$.
By Theorem~\ref{T:basic}, it follows that $A_p(w)<\infty$
and that $\ell^p(\ZZ,w)$ is an algebra.

\subsection{The finite-replacement lemma}

The idea of the proof of Theorem~\ref{T:main}
is to show that, starting with a weight $w$
such that $w=w_0$ outside a finite set, 
we can modify $w$ on a finite set to obtain a new
weight $\tilde{w}$ such that  $B_p(\tilde{w})$
is large, yet $A_p(\tilde{w})$
is hardly any bigger than $A_p(w)$.
Since  also $\tilde{w}=w_0$ outside a finite set, 
this process can be iterated, leading eventually  to a weight
satisfying the conclusions of Theorem~\ref{T:main}.

Here is the precise statement of the key step in this process.
We call it the finite-replacement lemma.

\begin{lemma}\label{L:perturb}
Let $w$ be a weight on $\ZZ$ such that $w=w_0$
outside a finite set. 
Given $\epsilon>0, M>0$ and $R>0$,
there exists a weight $\tilde{w}$ on~$\ZZ$ such that:
\begin{enumerate}[\normalfont(i)]
\item $\tilde{w}=w$ outside a finite subset of $(R,\infty)$;
\item $A_p(\tilde{w})\le A_p(w)+\epsilon$;
\item $\|1/\tilde{w}\|_q\le \|1/w\|_q+\epsilon$;
\item $B_p(\tilde{w})>M$.
\end{enumerate}
\end{lemma}

\subsection{Deduction of the main result}

Before going any further, let us see why this lemma implies
the theorem.

\begin{proof}[Proof of Theorem~\ref{T:main} (assuming Lemma~\ref{L:perturb})]
Choose sequences of weights $(w_n)$ 
and positive integers $(R_n)$
recursively as follows. 
Set $w_0(x):=(1+|x|)^s$ and $R_0:=1$.
Suppose that $n\ge1$ and that $w_0,\dots,w_{n-1}$ and $R_0,\dots,R_{n-1}$ have already been chosen.
Using Lemma~\ref{L:perturb}, choose $w_n$ such that:
\begin{itemize}
\item $w_n=w_{n-1}$ outside a finite subset of $(R_{n-1},\infty)$;
\item $A_p(w_{n})\le A_p(w_{n-1})+2^{-n}$;
\item $\|1/w_n\|_q\le\|1/w_{n-1}\|_q+2^{-n}$;
\item $B_p(w_{n})>n$.
\end{itemize}
Then choose $R_{n}$ large enough so that $R_n\ge 2R_{n-1}$
and 
\[
\Bigl(\sup_{|x|\le R_n/2}\sum_{|y|\le R_{n}/2}
\Bigl(\frac{w_n(x)}{w_n(y)w_n(x-y)}\Bigr)^q\Bigr)^{1/q}>\frac{B_p(w_n)}{2}.
\]
Since $R_n\to\infty$, it is clear that, for each $x\in\ZZ$, the 
sequence $w_n(x)$ is eventually constant. Set $w(x):=\lim_{n\to\infty}w_n(x)$. We shall show that this weight~$w$  satisfies the conclusions of Theorem~\ref{T:main}.

First we show that $B_p(w)=\infty$. For each $n\ge1$,
we have $w_n=w_{n+1}=w_{n+2}=\cdots$ on $[-R_n,R_n]$,
so $w=w_n$ on $[-R_n,R_n]$. Hence
\begin{align*}
B_p(w)
&=\Bigl(\sup_{x\in\ZZ}\sum_{y\in\ZZ}\Bigl(\frac{w(x)}{w(y)w(x-y)}\Bigr)^q\Bigr)^{1/q}\\
&\ge \Bigl(\sup_{|x|\le R_n/2}\sum_{|y|\le R_{n}/2}\Bigl(\frac{w(x)}{w(y)w(x-y)}\Bigr)^q\Bigr)^{1/q}\\
&= \Bigl(\sup_{|x|\le R_n/2}\sum_{|y|\le R_{n}/2}\Bigl(\frac{w_n(x)}{w_n(y)w_n(x-y)}\Bigr)^q\Bigr)^{1/q}\\
&\ge \frac{B_p(w_n)}{2}\ge \frac{n}{2}.
\end{align*}
As $n$ can be  arbitrarily large, we deduce that $B_p(w)=\infty$,
as claimed.

Next we show that $1/w\in\ell^q(\ZZ)$ (this is not
part of the statement of the theorem, but it is a necessary intermediate
step). By induction, we have 
$\|1/w_n\|_q\le\|1/w_0\|_q+1-2^{-n}$ for all $n\ge0$.
Thus,  as $w=w_n$ on $[-R_n,R_n]$,
\[
\Bigl(\sum_{|x|\le R_n}w(x)^{-q}\Bigr)^{1/q}
=\Bigl(\sum_{|x|\le R_n}w_n(x)^{-q}\Bigr)^{1/q}
\le \|1/w_n\|_q\le \|1/w_0\|_q+1.
\]
Letting $n\to\infty$, we deduce that $1/w\in\ell^q(\ZZ)$, as claimed.

Lastly we show that $A_p(w)<\infty$.
Let $f,g$ be finitely supported elements of $\ell^p(\ZZ,w)$.
Choose $n$ large enough so that $f,g$ and $f*g$ are all 
supported in $[-R_n,R_n]$, on which we have $w=w_n$.
Then
\[
\|f*g\|_{p,w}=\|f*g\|_{p,w_n}\le A_p(w_n)\|f\|_{p,w_n}\|g\|_{p,w_n}=A_p(w_n)\|f\|_{p,w}\|g\|_{p,w}.
\]
By induction, we have $A_p(w_n)\le A_p(w_0)+1-2^{-n}$
for all $n$. Hence
\begin{equation}\label{E:ineqconv}
\|f*g\|_{p,w}\le (A_p(w_0)+1)\|f\|_{p,w}\|g\|_{p,w}.
\end{equation}

It remains to extend the inequality \eqref{E:ineqconv} to those $f,g\in\ell^p(\ZZ,w)$
that are not finitely supported. For this, we need to use the
fact that $1/w\in\ell^q(\ZZ)$. By H\"older's inequality, if
$h\in\ell^p(\ZZ,w)$, then 
\[
\|h\|_1\le \|hw\|_p\|1/w\|_q=\|h\|_{p,w}\|1/w\|_q.
\]
Thus $\ell^p(\ZZ,w)$ embeds continuously into $\ell^1(\ZZ)$.
Given arbitrary functions $f,g\in\ell^p(\ZZ,w)$,
there exist finitely supported $f_k,g_k$ such that
$\|f_k-f\|_{p,w}\to0$ and $\|g_k-g\|_{p,w}\to0$.
By the preceding remarks, we have $f,g\in\ell^1(\ZZ)$
and $\|f_k-f\|_1\to0$ and $\|g_k-g\|_1\to0$.
Since $\ell^1(\ZZ)$ is a Banach algebra with respect 
to convolution, it follows that $f*g\in\ell^1(\ZZ)$
and $\|f_k*g_k-f*g\|_1\to0$.
In particular, $f*g$ is well-defined pointwise, and
$(f_k*g_k)\to (f*g)$ pointwise on $\ZZ$.
The inequality \eqref{E:ineqconv} now follows easily, using Fatou's lemma.
Thus $A_p(w)\le A_p(w_0)+1<\infty$, as claimed.
\end{proof}

All that now remains is to establish Lemma~\ref{L:perturb}.

%%%%%%%%%%%%%%%%%%%%%%%%%%%%%%%%%

\section{Proof of the finite-replacement lemma}\label{S:proof}

Throughout this section, we again
fix $1<p<\infty$, and set 
$q:=p/(p-1)$ and $s:=1/q+2/3$.
Note that $2/3<s<5/3$.
We write $w_0$ for the weight on $\ZZ$
defined by
$w_0(x):=(1+|x|)^s$.
We  suppose that, as in the lemma, we are given
$\epsilon>0, M>0, R>0$ and  a
weight $w$ on~$\ZZ$ such that $w=w_0$
outside a finite set. Thus, there exists a constant $C_w>0$
such that
\begin{equation}\label{E:w/w_0}
\frac{1}{C_w}\le \frac{w(x)}{w_0(x)}\le C_w \quad(x\in\ZZ).
\end{equation}

\subsection{Construction of $\tilde{w}$}

For convenience, we fix $L:= 32$;
this will act as a base integer
in our construction.
Let $N$ and $X$ be positive integers satisfying
\begin{equation}\label{E:initial}
N\ge 2 \quad\text{and}\quad X\ge 4000 N(2L)^N.
\end{equation}
The precise choice of $N$ and $X$ will be determined later.

Define $\phi:\ZZ\to[1,\infty)$  and $I\subset \ZZ$
by
\begin{equation}\label{E:phidef}
\phi(x):=\max\Bigl\{1,\,L^{N/2}e^{-|x|/L^N}\Bigr\} 
\quad\text{and}\quad
I:=\{x\in\ZZ:\phi(x)>1\}.
\end{equation}
A simple calculation shows that
$I\subset[-2NL^N,2NL^N]$. 

Now let $Y:=40NL^N$, and, for $1\le j\le N$, set 
\[
p_j:=X+Y2^j 
\quad\text{and}\quad
q_j:=2X-Y2^j.
\]
Thus $p_j+q_j=3X$ for all $j$. Put
\begin{equation}\label{E:Edef}
P_j:=p_j+I,\quad Q_j:=q_j+I
\quad\text{and}\quad
E:=\cup_{j=1}^N (P_j\cup Q_j).
\end{equation}
Note that, since $X\ge100Y\cdot 2^N$, 
the sets $P_j,Q_j~(1\le j \le N)$ are all
pairwise disjoint. 
The same inequality also implies that
$E\subset (X,2X)$.

Define $a:\ZZ\to[1,\infty)$ by
\[
a(x):=\begin{cases}
\phi(x-p_j), &x\in P_j,\,1\le j\le N,\\
\phi(x-q_j), &x\in Q_j,\,1\le j\le N,\\
1, &x\in\ZZ\setminus E.
\end{cases}
\]
Finally, we define  $\tilde{w}:\ZZ\to(0,\infty)$ by
\[
\tilde{w}:=w/a.
\]

In the next three subsections we shall estimate 
$B_p(\tilde{w}),\|1/\tilde{w}\|_q$ and $A_p(\tilde{w})$.
Using these estimates, we
shall then show that, with a suitable choice of $N$ and $X$,
the modified weight $\tilde{w}$ satisfies all the conditions (i)--(iv) of 
Lemma~\ref{L:perturb}. 

%%%%%%%%%%%%%%%%%%%%%%%%%%%%%%

\subsection{Estimate for $B_p(\tilde{w})$}

\begin{proposition}\label{P:Bpest}
\begin{equation}\label{E:Bpest}
B_p(\tilde{w})\ge \frac{L^{(1+1/q)N}N^{1/q}}{3^se^2C_w^3X^s}.
\end{equation}
\end{proposition}

\begin{proof}
Set $x_0:=3X$. Then $x_0\notin E$, so 
\[
a(x_0)=1.
\]
For $j=1,\dots,N$ and $1\le r\le L^N$, let  $y:=p_j+r=X+Y2^j+r$.
According to the definition of $a$, we have
\[
a(y)=\phi(r)=L^{N/2}e^{-|r|/L^N}\ge L^{N/2}/e.
\]
Also $x_0-y=q_j-r=2X-Y2^j-r$, so
\[
a(x_0-y)=\phi(-r)=L^{N/2}e^{-|r|/L^N}\ge L^{N/2}/e.
\]
Further, from \eqref{E:w/w_0}, we have
\begin{align*}
w(x_0)&\ge C_w^{-1}w_0(x_0)\ge C_w^{-1}(3X)^s,\\
w(y)&\le C_w w_0(y)\le C_w(3X)^s,\\
w(x_0-y)&\le C_w w_0(x_0-y)\le C_w(3X)^s.
\end{align*}
Hence, for each such $y$,
\[
\frac{\tilde{w}(x_0)}{\tilde{w}(x_0-y)\tilde{w}(y)}
=\frac{a(y)a(x_0-y)}{a(x_0)}\frac{{w}(x_0)}{{w}(x_0-y){w}(y)}
\ge\frac{L^N/e^2}{C_w^3(3X)^s}.
\]
There are $NL^N$  choices of such $y$, all of them distinct.
Hence
\[
\sum_{y\in\ZZ}\Bigl(\frac{\tilde{w}(x_0)}{\tilde{w}(x_0-y)\tilde{w}(y)}\Bigr)^q
\ge NL^N\frac{L^{Nq}}{C_w^{3q}(3X)^{sq}e^{2q}}.
\]
This implies that
\[
\Bigl(\sup_{x\in\ZZ}\sum_{y\in\ZZ}\Bigl(\frac{\tilde{w}(x)}{\tilde{w}(x-y)\tilde{w}(y)}\Bigr)^q\Bigr)^{1/q}
\ge  \frac{L^{(1+1/q)N}N^{1/q}}{C_w^3(3X)^se^2},
\]
which gives \eqref{E:Bpest}.
\end{proof}

%%%%%%%%%%%%%%%%%%%%%%%%%%%%%%%%

\subsection{Estimate for $\|1/\tilde{w}\|_q$}

\begin{proposition}\label{P:qnormest}
\begin{equation}\label{E:qnormest}
\|1/\tilde{w}\|_q\le \|1/w\|_q+8C_wX^{-s}N^{1/q}L^{(1/q+1/2)N}.
\end{equation}
\end{proposition}

To prove this result, we need 
an auxiliary lemma. 
Recall that $\phi$ and $I$ were defined in
\eqref{E:phidef}.

\begin{lemma}\label{L:phiineq}
\begin{equation}\label{E:phiineq}
\sum_{x\in I}\phi(x)^q\le 4L^{(1+q/2)N}.
\end{equation}
\end{lemma}

\begin{proof}
Summing the geometric progression, we have
\[
\sum_{x\in I}\phi(x)^q
\le 2\sum_{x=0}^\infty L^{qN/2}e^{-qx/L^N}
=\frac{2L^{qN/2}}{1-e^{-q/L^N}}.
\]
Since the function $t\mapsto t/(1-e^{-t})$ is increasing for $t>0$,
we have 
\[
\frac{1}{1-e^{-q/L^N}}=\frac{L^N}{q}\frac{q/L^N}{1-e^{-q/L^N}}
\le \frac{L^N}{q}\frac{q}{1-e^{-q}}=\frac{L^N}{1-e^{-q}}\le  2L^N.
\]
The inequality \eqref{E:phiineq} follows.
\end{proof}

\begin{corollary}\label{C:aineq}
\begin{equation}\label{E:aineq}
\sum_{x\in E}a(x)^q
\le 8N L^{(1+q/2)N}.
\end{equation}
\end{corollary}

\begin{proof}
From the definition of $E$, we have
\[
\sum_{x\in E}a(x)^q
=\sum_{j=1}^N\Bigl(\sum_{x\in P_j}a(x)^q+\sum_{x\in Q_j}a(x)^q\Bigr)
=2N\sum_{x\in I}\phi(x)^q.
\]
Combining this identity with \eqref{E:phiineq}
gives the result.
\end{proof}

\begin{proof}[Proof of Proposition~\ref{P:qnormest}]
As $\tilde{w}=w$ off $E$, the  triangle inequality gives
\[
\|1/\tilde{w}\|_q\le \|1/w\|_q+\Bigl(\sum_{x\in E}\tilde{w}(x)^{-q}\Bigr)^{1/q},
\]
and so it suffices to estimate this last sum.

For $x\in E$, we have 
$w(x)\ge  w_0(x)/C_w=(1+|x|)^s/C_w\ge X^{s}/C_w$.
Hence
\[
\sum_{x\in E} \tilde{w}(x)^{-q}
=\sum_{x\in E}\frac{w(x)^{-q}}{a(x)^{-q}}
\le C_w^qX^{-sq}\sum_{x\in E}a(x)^q.
\]
By Corollary~\ref{C:aineq}, we have
\[
\sum_{x\in E}a(x)^q
\le 8N L^{(1+q/2)N}.
\]
Hence
\[
\Bigl(\sum_{x\in E}\tilde{w}(x)^{-q}\Bigr)^{1/q}
\le C_wX^{-s}(8N)^{1/q}L^{(1/q+1/2)N},
\]
which implies the result.
\end{proof}

%%%%%%%%%%%%%%%%%%%%%%%%%%%%%%

\subsection{Estimate for $A_p(\tilde{w})$}

Our goal in this subsection is to establish the following result.

\begin{proposition}\label{P:Apest} Writing $\gamma:=\max\{0,\,(p-2)/p\}$,
we have
\begin{align}
A_p(\tilde{w})
&\le e^{1/L^{N/8}}A_p(w)\label{E:Apest}\\
&\quad+2^{7}C_w^3\Bigl(X^{-s}N^\gamma L^{(1+1/q)N}+X^{-s}N^{1/q}L^{(1/q+1/2)N}+L^{-N/12}\Bigr).\nonumber
\end{align}
\end{proposition}

This is the most complicated estimate.
The usual method for bounding $A_p(\tilde{w})$ via the
inequality $A_p(\tilde{w})\le B_p(\tilde{w})$ is of little help,
since our aim is to show that
$A_p(\tilde{w})$ is significantly smaller than $B_p(\tilde{w})$.
Instead we shall adopt a `divide-and-conquer' approach.
To facilitate this, we introduce some notation.

Let $v:\ZZ\to(0,\infty)$ be a weight such that $A_p(v)<\infty$.
Given non-empty subsets $S,T$ of $\ZZ$, we define
\begin{align*}
&A_p(v;S,T):=\\
&\qquad\sup\Bigl\{\|f*g\|_{p,v}: \|f\|_{p,v}\le 1,\,\|g\|_{p,v}\le1,\,\supp f\subset S,\,\supp g\subset T\Bigr\}.
\end{align*}
Clearly $A_p(v;\ZZ,\ZZ)=A_p(v)$.
Also, since convolution on $\ZZ$ is commutative,
\begin{equation}\label{E:symmetric}
A_p(v;S,T)=A_p(v;T,S) \quad(S,T\subset\ZZ).
\end{equation}
Further, by the triangle inequality for $\|\cdot\|_{p,v}$, we have
\begin{equation}\label{E:subadd}
\begin{aligned}
A_p(v;S_1\cup S_2,\,T)&\le A_p(v;S_1,T)+A_p(v;S_2,T) \quad(S_1,S_2,T\subset\ZZ),\\
A_p(v;S,\, T_1\cup T_2)&\le A_p(v;S, T_1)+A_p(v;S, T_2)
\quad(S,T_1,T_2\subset\ZZ).
\end{aligned}
\end{equation}

It will be convenient to express $A_p(w;S,T)$
in terms of unweighted $\ell^p$-norms.
Among other things, this will facilitate the use 
of the H\"older and Young inequalities \eqref{E:Holder}
and \eqref{E:Young} in
the arguments that follow.

\begin{lemma}\label{L:operator}
Let $v$ be a weight on $\ZZ$ such that $A_p(v)<\infty$.
Define $K_v:\ZZ\times\ZZ\to(0,\infty)$ by
\begin{equation}\label{E:Kdef}
K_{v}(y,z):=\frac{v(y+z)}{v(y)v(z)}
\quad(y,z\in\ZZ),
\end{equation}
and $\Phi_v:\ell^p(\ZZ)\times\ell^p(\ZZ)\to\ell^p(\ZZ)$ by
\begin{equation}\label{E:Phidef}
\Phi_v(F,G)(x):=\sum_{y+z=x}K_v(y,z)F(y)G(z)
\quad(F,G\in\ell^p(\ZZ)).
\end{equation}
Then $\Phi_v$ is a bounded bilinear map, and,
for each pair of subsets $S,T\subset\ZZ$,
\begin{equation}\label{E:operator}
A_p(v;S,T)=\bigl\|\Phi_v:\ell^p(S)\times\ell^p(T)\to\ell^p(\ZZ)\bigr\|,
\end{equation}
the operator norm of the restriction of $\Phi_v$
to $\ell^p(S)\times\ell^p(T)$.
\end{lemma}

\begin{proof}
Given $f,g\in\ell^p(\ZZ,v)$,  set
 $F:=fv$ and $G:=gv$.
Then
$\|f\|_{p,v}=\|F\|_p$ and $\|g\|_{p,v}=\|G\|_p$.
Also
\[
\|f*g\|_{p,v}=\Bigl\|v\cdot\Bigl(\frac{F}{v}*\frac{G}{v}\Bigr)\Bigr\|_p
=\bigl\|\Phi_v(F,G)\bigr\|_p.
\]
The result follows easily from these observations.
\end{proof}

We now undertake the estimation of $A_p(\tilde{w})$.
Our first result encapsulates the `divide-and-conquer' 
strategy mentioned earlier.

\begin{lemma}\label{L:dandc}
Let $H$ be an integer with $0<H\le L^N$.
Set $E^c:=\ZZ\setminus E$ and $E^c_H:=E^c\setminus [-H,H]$.
Then
\begin{equation}\label{E:dandc}
A_p(\tilde{w})\le e^{H/L^N}A_p(w)+A_p(\tilde{w}; E,E)+
2A_p(\tilde{w}; E,E_H^c).
\end{equation}
\end{lemma}

\begin{proof}
We begin with a preliminary remark.
We claim that, defining $K_w$ and $K_{\tilde{w}}$
as in \eqref{E:Kdef}, we have
\begin{equation}\label{E:claim}
K_{\tilde{w}}\le e^{H/L^N}K_w 
\quad\text{on~} \bigl(E^c\times E^c\bigr)\cup\bigl(E\times[-H,H]\bigr)\cup\bigl([-H,H]\times E\bigr).
\end{equation}
Indeed, from the definition of $K_w,K_{\tilde{w}}$, we have
\[
\frac{K_{\tilde{w}}(y,z)}{K_w(y,z)}=\frac{a(y)a(z)}{a(y+z)}
\qquad(y,z\in\ZZ).
\]
Since $a=1$ on $E^c$
and $a\ge 1$ everywhere, 
it follows that $K_{\tilde{w}}\le K_w$
on $E^c\times E^c$. Further, if $y\in E$ and $z\in[-H,H]$,
then $a(y)/a(y+z)\le e^{H/L^N}$ and $a(z)=1$, so $K_{\tilde{w}}\le e^{H/L^N}K_w$ on $E\times [-H,H]$. Likewise $K_{\tilde{w}}\le e^{H/L^N}K_w$ on $[-H,H]\times E$. This justifies the claim \eqref{E:claim}.

Now let $F,G\in\ell^p(\ZZ)$,
and let $\Phi_{\tilde{w}}$ be defined as in \eqref{E:Phidef}. 
Since $\ZZ\times \ZZ$ can be partitioned as
\[
\bigl(E\times E\bigr)\cup \bigl(E^c\times E^c\bigr)\cup\bigl(E\times [-H,H]\bigr)\cup\bigl([-H,H]\times E\bigr)\cup\bigl(E\times E^c_H\bigr)\cup\bigl(E^c_H\times E\bigr),
\]
the bilinearity of $\Phi_{\tilde{w}}$ implies that
\begin{equation}\label{E:longeqn}
\begin{aligned}
\Phi_{\tilde{w}}(F,G)
&=\Phi_{\tilde{w}}(F1_E,G1_E)+\Phi_{\tilde{w}}(F1_{E^c},G1_{E^c})+\Phi_{\tilde{w}}(F1_{E},G1_{[-H,H]})\\
&\quad+\Phi_{\tilde{w}}(F1_{[-H,H]},G1_E)+\Phi_{\tilde{w}}(F1_E,G1_{E^c_H})+
\Phi_{\tilde{w}}(F1_{E^c_H},G1_E).
\end{aligned}
\end{equation}
For each $x\in\ZZ$,
\begin{align*}
&\Bigl|\Phi_{\tilde{w}}(F1_{E^c},G1_{E^c})(x)+
\Phi_{\tilde{w}}(F1_{E},G1_{[-H,H]})(x)+
\Phi_{\tilde{w}}(F1_{[-H,H]},G1_{E})(x)\Bigr|\\
&=\Bigl|\sum_{y+z=x}K_{\tilde{w}}(y,z)F(y)G(z)
\Bigl(1_{E^c\times E^c}+1_{E\times[-H,H]}+1_{[-H,H]\times E}\Bigr)(y,z)\Bigr|\\
&\le \sum_{y+z=x}K_{\tilde{w}}(y,z)|F(y)||G(z)|
\Bigl(1_{E^c\times E^c}+1_{E\times[-H,H]}+1_{[-H,H]\times E}\Bigr)(y,z)\\
&\le \sum_{y+z=x}e^{H/L^N}K_w(y,z)|F(y)||G(z)| \qquad\text{(by \eqref{E:claim})}\\
&=e^{H/L^N}\Phi_w(|F|,|G|)(x).
\end{align*}
Inserting this information into \eqref{E:longeqn}, we deduce that
\begin{align*}
|\Phi_{\tilde{w}}(F,G)|&\le e^{H/L^N}\Phi_w(|F|,|G|)+|\Phi_{\tilde{w}}(F1_E,G1_E)|\\
&\qquad+|\Phi_{\tilde{w}}(F1_E,G1_{E^c_H})|+|\Phi_{\tilde{w}}(F1_{E^c_H},G1_E)|.
\end{align*}
Taking $\|\cdot\|_p$-norms and then using \eqref{E:operator},
we obtain
\begin{align*}
\|\Phi_{\tilde{w}}(F,G)\|_p&\le e^{H/L^N}\|\Phi_w(|F|,|G|)\|_p+\|\Phi_{\tilde{w}}(F1_E,G1_E)\|_p\\
&\qquad+\|\Phi_{\tilde{w}}(F1_E,G1_{E^c_H})\|_p+\|\Phi_{\tilde{w}}(F1_{E^c_H},G1_E)\|_p\\
&\le e^{H/L^N}A_p(w)\|F\|_p\|G\|_p
+A_p(\tilde{w};E,E)\|F\|_p\|G\|_p\\
&\qquad+A_p(\tilde{w};E,E^c_H)\|F\|_p\|G\|_p
+A_p(\tilde{w};E^c_H,E)\|F\|_p\|G\|_p\\
&=\Bigl(e^{H/L^N}A_p(w)+A_p(\tilde{w};E,E)+2A_p(\tilde{w};E,E^c_H)\Bigr)\|F\|_p\|G\|_p,
\end{align*}
where, in the last line, we used the symmetry
property \eqref{E:symmetric}. Thus
\[
\bigl\|\Phi_{\tilde{w}}:\ell^p(\ZZ)\times\ell^p(\ZZ)\to\ell^p(\ZZ)\bigr\|
\le e^{H/L^N}A_p(w)+A_p(\tilde{w};E,E)+2A_p(\tilde{w};E,E^c_H),
\]
which, in view of \eqref{E:operator}, gives \eqref{E:dandc}.
\end{proof}

It remains to estimate $A_p(\tilde{w};E,E)$ 
and $A_p(\tilde{w};E,E^c_H)$.
We begin with the easier one, namely $A_p(\tilde{w};E,E^c_H)$.

\begin{lemma}\label{L:EEch}
Let $H, E^c,E^c_H$ be as in Lemma~\ref{L:dandc}.
Then
\begin{equation}\label{E:EEch}
A_p(\tilde{w};E,E^c_H)
\le 32C_w^3\Bigl(X^{-s}N^{1/q}L^{(1/q+1/2)N}+L^{N/2}H^{-2/3}\Bigr).
 \end{equation}
\end{lemma}

\begin{proof}
Let $F,G\in\ell^p(\ZZ)$ 
with $\supp F\subset E$ and $\supp G\subset E^c_H$.
Defining $\Phi_{\tilde{w}}$ and $K_{\tilde{w}}$
as in Lemma~\ref{L:operator},
for each $x\in\ZZ$, we have
\begin{align*}
|\Phi_{\tilde{w}}(F,G)(x)|
&\le \sum_{y+z=x}K_{\tilde{w}}(y,z)|F(y)||G(z)|\\
&\le \sum_{y+z=x}\frac{a(y)a(z)}{a(y+z)}K_w(y,z)|F(y)||G(z)|\\
&\le \sum_{y+z=x}\frac{a(y)a(z)}{a(y+z)}C_w^3K_{w_0}(y,z)|F(y)||G(z)|\\
&\le \sum_{y+z=x}C_w^3K_{w_0}(y,z)|a(y)F(y)||G(z)|,
\end{align*}
the last inequality because $a(z)=1$ for all $z\in\supp G$ (which
is disjoint from $E$) and $a(y+z)\ge1$ for all $y,z$. Now
\begin{align*}
K_{w_0}(y,z)
&=\Bigl(\frac{1+|y+z|}{(1+|y|)(1+|z|)}\Bigr)^s
\le \Bigl(\frac{1}{1+|y|}+\frac{1}{1+|z|}\Bigr)^s\\
&\le 2^s\Bigl( \frac{1}{(1+|y|)^s}+\frac{1}{(1+|z|)^s}\Bigr)
=2^s\Bigl(\frac{1}{w_0(y)}+\frac{1}{w_0(z)}\Bigr).
\end{align*}
If $y\in \supp F$, then $y\in E$, so $|y|>X$, and hence $w_0(y)\ge X^s$.
Hence
\begin{align*}
|\Phi_{\tilde{w}}(F,G)(x)|
&\le 2^sC_w^3\sum_{y+z=x}\Bigl(\frac{1}{X^s}+\frac{1}{w_0(z)}\Bigr)|a(y)F(y)||G(z)|\\
&=2^sC_w^3\Bigl(\frac{1}{X^s}(|aF|*|G|)(x)+(|aF|*|G/w_0|)(x)\Bigr).
\end{align*}
Now, by Young's inequality and then H\"older's inequality,
\[
\Bigl\||aF|*|G|\Bigr\|_p\le \|aF\|_1\|G\|_p\le \|a1_E\|_q\|F\|_p\|G\|_p,
\]
where, from Corollary~\ref{C:aineq},
$\|a1_E\|_q\le (8N)^{1/q}L^{(1/q+1/2)N}$.
Also, by Young and H\"older again,
\[
\Bigl\||aF|*|G/w_0|\Bigr\|_p\le \|aF\|_p\|G/w_0\|_1
\le \|a\|_\infty\|(1/w_0)1_{\ZZ\setminus[-H,H]}\|_q\|F\|_p\|G\|_p,
\]
where, from the definition of $a$ and $\phi$, we have  $\|a\|_\infty\le L^{N/2}$, and where
\[
\|(1/w_0)1_{\ZZ\setminus[-H,H]}\|_q
=\Bigl(2\sum_{k>H}\frac{1}{(1+k)^{sq}}\Bigr)^{1/q}
\le 3^{1/q}H^{1/q-s}\le 3H^{-2/3}.
\]
Putting everything together, we obtain
\begin{align*}
\|\Phi_{\tilde{w}}(F,G)\|_p
&\le 2^sC_w^3\Bigl(X^{-s}(8N)^{1/q}L^{(1/q+1/2)N}+L^{N/2}3H^{-2/3}\Bigr)\|F\|_p\|G\|_p\\
&\le 32C_w^3\Bigl(X^{-s}N^{1/q}L^{(1/q+1/2)N}+L^{N/2}H^{-2/3}\Bigr)\|F\|_p\|G\|_p.
\end{align*}
In view of \eqref{E:operator}, this gives \eqref{E:EEch}.
\end{proof}

We now turn to the estimation of $A_p(\tilde{w};E,E)$.
This depends heavily on the internal structure of $E$.
Recall from \eqref{E:Edef} that $E=\cup_{j=1}^N(P_j\cup Q_j)$,
where $P_j,Q_j$ are particular translates of the interval
$I$ defined in \eqref{E:phidef}.
We first treat these translates individually.

\begin{lemma}\label{L:individual}
Let $S,T\in\{P_1,\dots,P_N,Q_1,\dots,Q_N\}$.
Then
\begin{equation}\label{E:individual}
A_p(\tilde{w};S,T)\le 16C_w^3X^{-s}L^{(1+1/q)N}.
\end{equation}
\end{lemma}

\begin{proof}
Let $y\in S$ and $z\in T$. Then, defining $K_{\tilde{w}}$ as in Lemma~\ref{L:operator}, we have
\[
K_{\tilde{w}}(y,z)
=\frac{\tilde{w}(y+z)}{\tilde{w}(y)\tilde{w}(z)}
=\frac{a(y)a(z)}{a(y+z)}\frac{w(y+z)}{w(y)w(z)}.
\]
Now $1/C_w\le w/w_0\le C_w$. Also, $y,z\in E\subset (X,2X)$.
It follows that
\[
\frac{w(y+z)}{w(y)w(z)}\le C_w^3\frac{w_0(y+z)}{w_0(y)w_0(z)}
\le 4C_w^3X^{-s}.
\]
Furthermore, $a(y+z)\ge1$. Therefore
\[
K_{\tilde{w}}(y,z)
\le 4C_w^3X^{-s}a(y)a(z)
\quad(y\in S,\, z\in T).
\]

Now let $F,G\in\ell^p(\ZZ)$ with $\supp F\subset S$
and $\supp G\subset T$. 
The above estimate for $K_{\tilde{w}}$
implies that, defining $\Phi_{\tilde{w}}$
as in Lemma~\ref{L:operator}, we have
\[
\|\Phi_{\tilde{w}}(F,G)\|_p\le 4C_w^3X^{-s}\bigl\||aF|*|aG|\bigr\|_p.
\]
By Young's inequality,
\[
\bigl\||aF|*|aG|\bigr\|_p\le \|aF\|_1\|aG\|_p.
\]
By H\"older's inequality, the definition of $a$, and Lemma~\ref{L:phiineq},
\begin{align*}
\|aF\|_1&\le \|a1_S\|_q\|F\|_p=\Bigl(\sum_{x\in S}|a(x)|^q\Bigr)^{1/q}\|F\|_p\\
&=\Bigl(\sum_{x\in I}|\phi(x)|^q\Bigr)^{1/q}\|F\|_p\le 4^{1/q}L^{(1/q+1/2)N}\|F\|_p.
\end{align*}
Also, clearly,
\[
\|aG\|_p\le \|a\|_\infty\|G\|_p=\|\phi\|_\infty\|G\|_p=L^{N/2}\|G\|_p.
\]
Combining these estimates gives
\[
\|\Phi_{\tilde{w}}(F,G)\|_p\le 4C_w^3X^{-s}4^{1/q}L^{(1+1/q)N}\|F\|_p\|G\|_p.
\]
Hence, using Lemma~\ref{L:operator}, we deduce that
\[
A_p(\tilde{w};S,T)=\bigl\|\Phi_{\tilde{w}}:\ell^p(S)\times\ell^p(T)\to\ell^p(\ZZ)\bigr\|\le 4C_w^3X^{-s}4^{1/q}L^{(1+1/q)N},
\]
whence the result.
\end{proof}

We now need to combine the individual estimates in
Lemma~\ref{L:individual} to obtain a bound for $A_p(\tilde{w};E,E)$.
One way to achieve this is to exploit the subadditivity property
\eqref{E:subadd} repeatedly to obtain
\begin{align*}
A_p(\tilde{w};E,E)
&\le \sum_{S,T\in\{P_1,\dots,P_N,Q_1,\dots,Q_N\}}A_p(\tilde{w};S,T)\\
&\le (2N)^2 16C_w^3X^{-s}L^{(1+1/q)N}.
\end{align*}
Although this is perfectly legitimate, it does not quite cut the mustard.
Indeed, it leads to an upper bound for $A_p(\tilde{w})$
that is bigger than a constant multiple of $X^{-s}N^2L^{(1+1/q)N}$,
whereas the lower bound  for $B_p(\tilde{w})$
already obtained in \eqref{E:Bpest} is a constant multiple of
$X^{-s}N^{1/q}L^{(1+1/q)N}$,
which is smaller by a factor of $N^{2-1/q}$. Thus, with these estimates alone, 
there is no hope of achieving our aim of showing that $A_p(\tilde{w})$ is
significantly smaller than $B_p(\tilde{w})$.

Instead, we need to make a more refined estimate of $A_p(\tilde{w};E,E)$ that takes advantage of the particular choice
of the blocks $P_j,Q_j$. 
It was already remarked, when these were
defined in \eqref{E:Edef}, that they are disjoint. In fact, they enjoy
the following stronger property.

\begin{lemma}\label{L:P+Q}
Let $i,j,k,\ell\in\{1,\dots,N\}$. Then
\begin{align*}
(P_i+P_j)\cap (P_k+P_\ell)\ne\emptyset &\Longrightarrow \{i,j\}=\{k,\ell\}\\
(Q_i+Q_j)\cap (Q_k+Q_\ell)\ne\emptyset &\Longrightarrow \{i,j\}=\{k,\ell\}\\
(P_i+Q_j)\cap (P_k+Q_\ell)\ne\emptyset &\Longrightarrow \{i,\ell\}=\{j,k\}.
\end{align*}
\end{lemma}

\begin{proof}
Recall from \eqref{E:Edef} that $P_i=p_i+I=X+Y2^i+I$ for
each $i$. Thus 
\begin{align*}
&(P_i+P_j)\cap (P_k+P_\ell)\ne\emptyset\\
&\qquad\Rightarrow (2X+Y(2^i+2^j)+2I)\cap(2X+Y(2^k+2^\ell)+2I)\ne\emptyset\\
&\qquad\Rightarrow (Y(2^i+2^j)+2I)\cap(Y(2^k+2^\ell)+2I)\ne\emptyset\\
&\qquad \Rightarrow Y(2^i+2^j-2^k-2^\ell)\in 4I.
\end{align*}
Now, as remarked after the definition \eqref{E:phidef} of $I$,
we have $I\subset[-2NL^N,\,2NL^N]$.
Also, $Y=40NL^N$, by definition.
Therefore, the only way that we can have 
$Y(2^i+2^j-2^k-2^\ell)\in 4I$ is if $2^i+2^j=2^k+2^\ell$.
By uniqueness of binary expansion, this implies that
$\{i,j\}=\{k,\ell\}$. This establishes the first implication  in
the lemma.

The second implication is very similar, and we omit the proof.

For the third implication, we need the definitions of both
$P_i$ and $Q_j$, namely $P_i=p_i+I=X+Y2^i+I$
and $Q_j=2X-Y2^j+I$. Thus
\begin{align*}
&(P_i+Q_j)\cap (P_k+Q_\ell)\ne\emptyset\\
&\qquad\Rightarrow (3X+Y(2^i-2^j)+2I)\cap(3X+Y(2^k-2^\ell)+2I)\ne\emptyset\\
&\qquad\Rightarrow (Y(2^i-2^j)+2I)\cap(Y(2^k-2^\ell)+2I)\ne\emptyset\\
&\qquad \Rightarrow Y(2^i-2^j-2^k+2^\ell)\in 4I.
\end{align*}
Just as before, the containment $Y(2^i-2^j-2^k+2^\ell)\in 4I$
can only hold if $2^i+2^\ell=2^j+2^k$,
which implies that $\{i,\ell\}=\{j,k\}$. 
This proves the third implication in the lemma.
\end{proof}

We shall also need the following elementary inequality.

\begin{lemma}\label{L:elem}
Let $\alpha_1,\dots,\alpha_N,\beta_1,\dots,\beta_N\in[0,\infty)$.
Then
\[
\sum_{i=1}^N\alpha_i\beta_i
\le N^\gamma\Bigl(\sum_{i=1}^N\alpha_i^p\Bigr)^{1/p}\Bigl(\sum_{i=1}^N\beta_i^p\Bigr)^{1/p},
\]
where $\gamma:=\max\bigl\{0,\,(p-2)/p\bigr\}$.
\end{lemma}

\begin{proof}
If $p\le 2$, then this follows by applying H\"older's inequality
to $\sum_{i}\alpha_i\beta_i$, and then using the fact that,
since $p\le q$, we have
$\|\cdot\|_q\le \|\cdot\|_p$ on $\RR^N$.

If $p>2$, then it follows directly by applying H\"older to $\sum_i(\alpha_i\cdot\beta_i\cdot1)$
with indices $p,p,p/(p-2)$.
\end{proof}

We can now establish the refined decomposition estimate
for $A_p(\tilde{w};E,E)$.

\begin{lemma}\label{L:decomp}
\begin{align*}
A_p(\tilde{w};E,E)
&\le 
2\max_{i,j}A_p(\tilde{w};P_i,P_j)
+2\max_{i,j}A_p(\tilde{w};Q_i,Q_j)\\
&\qquad+4N^\gamma\max_{i,j}A_p(\tilde{w};P_i,Q_j).
\end{align*}
\end{lemma}

\begin{proof}
Setting $P:=\cup_{j=1}^NP_j$ and $Q:=\cup_{j=1}^NQ_j$,
we have $E=P\cup Q$, so, by the subadditivity property
\eqref{E:subadd},
\[
A_p(\tilde{w};E,E)\le A_p(\tilde{w};P,P)+A_p(\tilde{w};Q,Q)+A_p(\tilde{w};P,Q)+A_p(\tilde{w};Q,P).
\]
We shall estimate separately each of the quantities on the right-hand side.

Let us begin with  $A_p(\tilde{w};P,P)$.
Let $f,g\in\ell^p(\ZZ,\tilde{w})$ with 
$\supp f\subset P $ and $\supp g\subset P$.
Then
\[
f*g=\sum_{i,j=1}^N(f1_{P_i})*(g1_{P_j})
=\sum_{1\le i<j\le N}h_{i,j}+\sum_{1\le i\le N}h_i,
\]
where 
\[
h_{i,j}:=(f1_{P_i})*(g1_{P_j})+(f1_{P_j})*(g1_{P_i})
\quad\text{and}\quad
h_i:=(f1_{P_i})*(g1_{P_i}).
\]
The support of $h_{i,j}$ is contained in $P_i+P_j$,
and that of $h_i$ in $P_i+P_i$.
Thus, by Lemma~\ref{L:P+Q}, the functions 
$h_{i,j}$ and $h_i$ all have disjoint supports.
Therefore
\begin{align*}
|f*g|^p
&\le \sum_{1\le i<j\le N}|h_{i,j}|^p+\sum_{1\le i\le N}|h_i|^p\\
&\le \sum_{1\le i<j\le N}2^p\Bigl(|(f1_{P_i})*(g1_{P_j})|^p+|(f1_{P_j})*(g1_{P_i})|^p\Bigr)\\
&\qquad+\sum_{1\le i\le N}|(f1_{P_i})*(g1_{P_i})|^p\\
&\le 2^p\sum_{i,j=1}^N|(f1_{P_i})*(g1_{P_j})|^p.
\end{align*}
It follows that
\begin{align*}
\|f*g\|_{p,\tilde{w}}^p
&\le 2^p\sum_{i,j=1}^N\|(f1_{P_i})*(g1_{P_j})\|_{p,\tilde{w}}^p\\ 
&\le 2^p\sum_{i,j=1}^N A_p(\tilde{w};P_i,P_j)^p\|f1_{P_i}\|_{p,\tilde{w}}^p\|g1_{P_j}\|_{p,\tilde{w}}^p\\
&\le 2^p\max_{i,j}A_p(\tilde{w};P_i,P_j)^p\sum_{i,j=1}^N \|f1_{P_i}\|_{p,\tilde{w}}^p\|g1_{P_j}\|_{p,\tilde{w}}^p\\
&= 2^p\max_{i,j}A_p(\tilde{w};P_i,P_j)^p\|f\|_{p,\tilde{w}}^p\|g\|_{p,\tilde{w}}^p,
\end{align*}
where, in the last equality, we used the fact that
the sets $P_1,\dots,P_N$ are pairwise disjoint.
Hence
\[
A_p(\tilde{w};P,P)\le 2\max_{i,j}A_p(\tilde{w};P_i,P_j).
\]

Likewise we have
\[
A_p(\tilde{w};Q,Q)\le  2\max_{i,j}A_p(\tilde{w};Q_i,Q_j).
\]
The proof is almost identical, so we omit the details.

The case of $A_p(\tilde{w};P,Q)$ is a little  different,
due to the fact that 
the sets $P_j+Q_j~(j=1,\dots,N)$ all coincide.

Let $f,g\in\ell^p(\ZZ,\tilde{w})$ with 
$\supp f\subset P $ and $\supp g\subset Q$.
Then
\[
f*g=\sum_{i,j=1}^N(f1_{P_i})*(g1_{Q_j})
=\sum_{\substack{i,j=1\\i\ne j}}^Nk_{i,j}+k,
\]
where 
\[
k_{i,j}:=(f1_{P_i})*(g1_{Q_j})
\quad\text{and}\quad
k:=\sum_{i=1}^N(f1_{P_i})*(g1_{Q_i}).
\]
The support of $k_{i,j}$ is contained in $P_i+Q_j$,
and that of $k$ in $P_1+Q_1$.
Thus, by Lemma~\ref{L:P+Q}, the functions 
$k_{i,j}$ and $k$ all have disjoint supports.
Therefore
\[
|f*g|^p
\le \sum_{\substack{i,j=1\\i\ne j}}^N |k_{i,j}|^p+|k|^p
= \sum_{\substack{i,j=1\\i\ne j}}^N \bigl|(f1_{P_i})*(g1_{Q_j})\bigr|^p+\Bigl|\sum_{i=1}^N(f1_{P_i})*(g1_{Q_i})\Bigr|^p.
\]
It follows that
\[
\|f*g\|_{p,\tilde{w}}^p\le
\sum_{\substack{i,j=1\\i\ne j}}^N \bigl\|(f1_{P_i})*(g1_{Q_j})\bigr\|_{p,\tilde{w}}^p+
\Bigl\|\sum_{i=1}^N(f1_{P_i})*(g1_{Q_i})\Bigr\|_{p,\tilde{w}}^p.
\]
The first term on the right-hand side is treated
in the same way as in the previous case:
\begin{align*}
\sum_{\substack{i,j=1\\i\ne j}}^N\bigl\|(f1_{P_i})*(g1_{Q_j})\bigr\|_{p,\tilde{w}}^p
&\le \sum_{\substack{i,j=1\\i\ne j}}^N A_p(\tilde{w};P_i,Q_j)^p\|f1_{P_i}\|_{p,\tilde{w}}^p\|g1_{Q_j}\|_{p,\tilde{w}}^p\\
&\le  \max_{\substack{i,j\\i\ne j}}A_p(\tilde{w};P_i,Q_j)^p\|f\|_{p,\tilde{w}}^p\|g\|_{p,\tilde{w}}^p.
\end{align*}
However, the second term has to be treated differently,
since the sets $P_i+Q_i$ are all the same.
Using the triangle inequality for $\|\cdot\|_{p,\tilde{w}}$,
the definition of $A_p(\tilde{w};P_i,Q_i)$,
and then Lemma~\ref{L:elem}, we have
\begin{align*}
&\Bigl\|\sum_{i=1}^N(f1_{P_i})*(g1_{Q_i})\Bigr\|_{p,\tilde{w}}^p\\
&\qquad\le \Bigl(\sum_{i=1}^N\bigl\|(f1_{P_i})*(g1_{Q_i})\bigr\|_{p,\tilde{w}}\Bigr)^p\\
&\qquad\le \Bigl(\sum_{i=1}^N A_p(\tilde{w};P_i,Q_i)\|f1_{P_i}\|_{p,\tilde{w}}\|g1_{Q_i}\|_{p,\tilde{w}}\Bigr)^p\\
&\qquad\le \max_i A_p(\tilde{w};P_i,Q_i)^p \Bigl(\sum_{i=1}^N \|f1_{P_i}\|_{p,\tilde{w}}\|g1_{Q_i}\|_{p,\tilde{w}}\Bigr)^p\\
&\qquad\le \max_i A_p(\tilde{w};P_i,Q_i)^p N^{p\gamma}\Bigl(\sum_{i=1}^N \|f1_{P_i}\|_{p,\tilde{w}}^p\Bigr)\Bigl(\sum_{i=1}^N\|g1_{Q_i}\|_{p,\tilde{w}}^p\Bigr)\\
&\qquad=\max_i A_p(\tilde{w};P_i,Q_i)^p N^{p\gamma}\|f\|_{p,\tilde{w}}^p\|g\|_{p,\tilde{w}}^p.
\end{align*}
Putting everything together, we arrive at the estimate
\[
\|f*g\|_{p,\tilde{w}}^p\le (1+N^{p\gamma})\max_{i,j} A_p(\tilde{w};P_i,Q_j)^p\|f\|_{p,\tilde{w}}^p\|g\|_{p,\tilde{w}}^p,
\]
whence it follows that
\[
A_p(\tilde{w};P,Q)\le 2N^\gamma \max_{i,j}A_p(\tilde{w};P_i,Q_j).
\]
By the symmetry property \eqref{E:symmetric},
\[
A_p(\tilde{w};Q,P)=A_p(\tilde{w};P,Q).
\]
This completes the proof.
\end{proof}

\begin{corollary}\label{C:EEest}
Let $\gamma:=\max\{0,\,(p-2)/p\}$.
Then
\begin{equation}\label{E:EEest}
A_p(\tilde{w};E,E)\le 
2^7C_w^3X^{-s}N^\gamma L^{(1+1/q)N}.
\end{equation}
\end{corollary}

\begin{proof}
Combine Lemmas~\ref{L:individual} and \ref{L:decomp}.
\end{proof}

At long last, we have all that we need for the desired bound for $A_p(\tilde{w})$.

\begin{proof}[Proof of Proposition~\ref{P:Apest}]
Combining Lemma~\ref{L:dandc}, Lemma~\ref{L:EEch} and Corollary~\ref{C:EEest},
we have
\begin{align*}
A_p(\tilde{w})&\le e^{H/L^N}A_p(w)+2^7C_w^3X^{-s}N^\gamma L^{(1+1/q)N}\\
&\qquad+2\cdot32C_w^3\Bigl(X^{-s}N^{1/q}L^{(1/q+1/2)N}+L^{N/2}H^{-2/3}\Bigr).
\end{align*}
Here we are free to choose $H$ to be  
any integer satisfying $0<H\le L^N$.
We would like to `balance' the terms $H/L^N$ and $L^{N/2}H^{-2/3}$,
so that both tend to zero as $N\to\infty$. 
A suitable choice is $H:=\lfloor L^{7N/8}\rfloor$.
Inserting this into the above inequality and tidying up, 
we arrive at \eqref{E:Apest}.
\end{proof}

%%%%%%%%%%%%%%%%%%%%%%%%%%%%%%

\subsection{Conclusion of the proof of Lemma~\ref{L:perturb}}

To establish Lemma~\ref{L:perturb}, we need to show that, 
with a suitable choice of $N$ and $X$ satisfying \eqref{E:initial},  
all the conditions~(i)--(iv)
of the lemma are  satisfied. 

Condition~(i) will be satisfied if, 
\[
X\ge R.
\]
Indeed,  we have $\tilde{w}=w$ outside $E$, which is a finite subset of $(X,2X)$. Thus, if  $X\ge R$,
then $\tilde{w}=w$ outside a finite
subset of $(R,\infty)$.

By Proposition~\ref{P:Apest}, condition~(ii) will be satisfied if
\[
(e^{1/L^{N/8}}-1)A_p(w)<\epsilon/2
\]
and
\[
2^{7}C_w^3\Bigl(X^{-s}N^\gamma L^{(1+1/q)N}+X^{-s}N^{1/q}L^{(1/q+1/2)N}+L^{-N/12}\Bigr)<\epsilon/2,
\]
where $\gamma:=\max\{0,\,(p-2)/p\}$.

By Proposition~\ref{P:qnormest}, condition~(iii) will be satisfied if
\[
8C_wX^{-s}N^{1/q}L^{(1/q+1/2)N}<\epsilon.
\]

Lastly, by Proposition~\ref{P:Bpest}, condition~(iv) will be satisfied if
\[
\frac{L^{(1+1/q)N}N^{1/q}}{3^se^2C_w^3X^s}>M.
\]

Thus, to complete the proof of Lemma~\ref{L:perturb},
it suffices to establish the following result.
As before, we write $\gamma:=\max\{0,\,(p-2)/p\}$.

\begin{lemma}\label{L:constants}
Given $\delta>0$ and $\Delta>0$, there exist positive integers
$N$ and $X$ such that
\begin{align}
X&\ge 4000N(2L)^N, \label{E:cond1}\\
X&\ge \Delta, \label{E:cond2}\\
e^{1/L^{N/8}}-1&\le \delta,\label{E:cond3}\\
L^{-N/12}&\le \delta,\label{E:cond3.5}\\
X^{-s}N^\gamma L^{(1+1/q)N}&\le \delta, \label{E:cond4}\\
X^{-s}N^{1/q}L^{(1/2+1/q)N}&\le \delta, \label{E:cond5}\\
X^{-s}N^{1/q}L^{(1+1/q)N}&\ge \Delta.\label{E:cond6}
\end{align}
\end{lemma}

\begin{proof}
We take $X=X_N$, where
\[
X_N:=\Bigl\lfloor \Bigl(\frac{N^{1/q}L^{(1+1/q)N}}{\Delta}\Bigr)^{1/s }\Bigr\rfloor.
\]
Thus \eqref{E:cond6} is automatically satisfied. 
We shall show that each of the other conditions \eqref{E:cond1}--\eqref{E:cond5} 
is also satisfied, provided that $N$ is sufficiently large.

First of all, since $X_N\to\infty$ as $N\to\infty$, it is clear that \eqref{E:cond2} holds for all $N$ large. Also,
as $(e^{1/L^{N/8}}-1)$ and $L^{-N/12}$
both tend to zero as $N\to\infty$, 
it is evident that \eqref{E:cond3} and \eqref{E:cond3.5} hold for all $N$ large.
Further, as $N\to\infty$, we have
\[
X_N^{-s}N^\gamma L^{(1+1/q)N}\sim N^{\gamma-1/q}\Delta
\to0
\]
and
\[
X_N^{-s}N^{1/q} L^{(1/2+1/q)N}\sim L^{-N/2}\Delta\to0,
\]
so conditions~\eqref{E:cond4} and \eqref{E:cond5}
both hold if $N$ is sufficiently large.

It remains to verify \eqref{E:cond1}. 
This is where the precise choices of $s$ and $L$ finally enter the
picture. We recall that $s=1/q+2/3$, so in particular $s<5/3$.
Also $L=32=2^5$. Hence, as $N\to\infty$,
\[
\frac{X_N}{N(2L)^N}\sim
\frac{N^{1/qs}L^{(1/s+1/qs)N}}{\Delta^{1/s}N (2L)^N}
\sim \Delta^{-1/s}N^{(1/qs-1)}2^{(5/s+5/qs-6)N},
\]
and the right-hand side tends to infinity with $N$ since
\[
\frac{5}{s}+\frac{5}{qs}-6=\frac{5}{3s}-1>0.
\]
It follows  that condition \eqref{E:cond1}
also holds for all $N$ sufficiently large.
\end{proof}

To complete the proof of Lemma~\ref{L:perturb}, 
choose $\delta>0$ so small that
\[
\delta A_p(w)<\frac{\varepsilon}{2},\qquad
3\cdot 2^7 C_w^3\delta<\frac{\varepsilon}{2},
\qquad
8C_w\delta<\varepsilon
\]
and 
\[
\Delta>\max\{R,\,3^se^2C_w^3M\}.
\]
Applying Lemma~\ref{L:constants} with these values of $\delta$ and $\Delta$ concludes the proof of Lemma~\ref{L:perturb},
and with it, that of the main result Theorem~\ref{T:main}.

%%%%%%%%%%%%%%%%%%%%%%%%%%%%%%

\section{Concluding remarks}\label{S:conclusion}

\subsection{A generalization to infinite, finitely generated abelian groups}

The main result Theorem~\ref{T:main} actually implies a stronger version of itself  that is worth recording.

\begin{theorem}\label{T:stronger}
Let $1<p<\infty$, and let $G$ be an infinite, finitely generated abelian group. Then there exists a weight $w:G\to(0,\infty)$ such that $\ell^p(G,w)$ is an algebra with respect to convolution, but,
with $q=p/(p-1)$, we have
\[
\sup_{x\in G}\sum_{y\in G}\Bigl(\frac{w(x)}{w(y)w(x-y)}\Bigr)^q=\infty.
\]
\end{theorem}

\begin{proof}
We use the following elementary product observation. If $G_1,G_2$ are discrete abelian groups, $w_i$ is a weight on $G_i~(i=1,2)$, and
\[
w(x_1,x_2):=w_1(x_1)w_2(x_2),
\]
then
\begin{equation}\label{E:prodids}
A_p(w)= A_p(w_1)A_p(w_2)
\quad\text{and}\quad
B_p(w)=B_p(w_1)B_p(w_2).
\end{equation}
The second identity follows by separating the sums in the definition of $B_p$. For the first, the lower bound follows by testing on tensor products. For the upper bound, apply the convolution estimate on $G_2$ to each pair of slices and then the convolution estimate on $G_1$ to the resulting $\ell^p(G_2)$-norms.

By the fundamental theorem of finitely generated abelian groups,
$G$ is isomorphic to $\ZZ^r\times F$,
where $r$ is a strictly positive integer and $F$ is a finite abelian group. Let $v$ be a weight on $\ZZ$ supplied by Theorem~\ref{T:main}, so that
$A_p(v)<\infty$ and
$B_p(v)=\infty$.
On $\ZZ^r\times F$, define
\[
w(x_1,\dots,x_r,y)
:=
\prod_{j=1}^r v(x_j)
\qquad (x_1,\dots,x_r\in\ZZ,\,y\in F).
\]
Thus $w$ is the product of $r$ copies of $v$ and the trivial weight on $F$. Since convolution on the finite group $F$ is bounded on $\ell^p(F)$, repeated use of the product relations \eqref{E:prodids} gives
\[
A_p(w)<\infty
\quad\text{and}\quad
B_p(w)=\infty.
\]
Transporting $w$ through the above group isomorphism proves the result.
\end{proof}

The preceding extension was already observed in \cite[Proposition~8.2]{Mo19} in the case $p=2$, conditional at the time on the existence of a counterexample on $\ZZ$. The present theorem therefore completes that argument and, together with the product observation above, extends it to every $1<p<\infty$.

%%%%%%%%%%%%%%%%%%%%%%%%%%%%%%%%%%

\subsection{How to characterize bounded convolution?}

Theorem~\ref{T:main}
shows that the condition $B_p(w)<\infty$ does not characterize the weights $w$ for which $\ell^p(\ZZ,w)$ is an algebra.
One might hope to obtain a complete characterization by replacing the 
$\ell^q$-norm in the definition of $B_p(w)$ by some other norm within the scale of $\ell^r$-norms. We shall show that this is not possible.

For $0<r\le\infty$, set
\[
B_{p,r}(w):=\sup_{x\in\ZZ} \left\|\frac{w(x)}{w(\cdot)\,w(x-\cdot)}\right\|_r,
\]
so that $B_p(w)=B_{p,q}(w)$.

If $r\le q$, then $B_{p,r}(w)\ge B_p(w)$, so the condition $B_{p,r}(w)<\infty$ is sufficient for $\ell^p(\ZZ,w)$ to be an algebra, but, by Theorem~\ref{T:main}, it is not necessary.

If $r>q$, then the condition $B_{p,r}(w)<\infty$ is not sufficient. Indeed, let $w(x):=(1+|x|)^t$, where $\frac{1}{r}<t\le \frac{1}{q}$. The calculation in \S\ref{S:prelim}, with $qs$ replaced by $rt$, shows that $B_{p,r}(w)<\infty$. 
On the other hand, $\ell^p(\ZZ,w)$ is not an algebra. To see this, suppose the contrary, and, for $n\ge1$, set $f_n:=w^{-q}1_{[-n,n]}$ and $S_n:=\sum_{|x|\le n}w(x)^{-q}$. 
Then $\|f_n\|_{p,w}=S_n^{1/p}$, 
and $f_n*f_n$ is a positive function supported in $[-2n,2n]$, 
so, by H\"older's inequality,
\[
S_n^2
=\|f_n*f_n\|_1 \le \|f_n*f_n\|_{p,w}\,S_{2n}^{1/q} 
\le A_p(w)\,S_n^{2/p}\,S_{2n}^{1/q}.
\]
As $w$ is even and increasing on $\ZZ^+$, 
we have $S_{2n}\le 2S_n$. 
It follows that $S_n\le 2A_p(w)^q$ for all $n$, 
which contradicts the fact that 
$\sum_{x\in\ZZ}w(x)^{-q}=\infty$ when $tq\le1$.

The reason for this failure to characterize $A_p(w)<\infty$
is that the gap between $A_p(\widetilde{w})$ and $B_p(\widetilde{w})$ obtained in \S\ref{S:proof} comes from the arithmetic structure of the set $E$ on which the weight is modified, as expressed in Lemma~\ref{L:P+Q}. A condition of the form $B_{p,r}(w)<\infty$ depends only on the sizes of the values taken by the functions $y\mapsto w(x)/(w(y)w(x-y))$, and not on where in $\ZZ$ these values are taken.

%%%%%%%%%%%%%%%%%%%%%%%%%%%%%%%%%

\section*{AI tools declaration}

In preparing this work, the authors used ChatGPT to develop and explore strategies for the proof of Theorem~\ref{T:main}. The final mathematical  text was written entirely by the authors, who thoroughly corrected, verified, and streamlined the proposed strategy. The authors accept full responsibility for the accuracy and integrity of the final results presented in this paper.

\bibliographystyle{plain}
\bibliography{biblist}

\end{document}